\documentclass{amsart}

\usepackage{amsmath,amssymb}
\usepackage{geometry}
\usepackage{graphicx}

\def\ep{\varepsilon}
\def\SSN{{\mathbb S}^{N-1}}

\def\RN{\mathbb{R}^N}

\def\cA{\mathcal{A}}

\def\cH{\mathcal{H}}

\def\al{\alpha}
\def\la{\lambda}
\def\cP{\mathcal P}
\def\cQ{\mathcal{Q}}

\def\pa{\partial}
\def\fhi{\varphi}

\def\dist{\mbox{\rm dist}}

\def\diam{\mbox{\rm diam}}

\newcommand{\core}{\heartsuit}
\newcommand{\cK}{\mathcal K}
\newcommand{\ccK}{\core(\cK)}
\newcommand{\De}{\Delta}
\newcommand{\cT}{\mathcal T}
\newcommand{\be}{\beta}
\newcommand{\cR}{\mathcal R}
\newcommand{\cM}{\mathcal M}
\newcommand{\ga}{\gamma}
\newcommand{\Om}{\Omega}
\newcommand{\om}{\omega}
\newcommand{\tht}{\vartheta}
\newcommand{\RE}{\mathbb R}

\newtheorem{lemma}{Lemma}[section]
\newtheorem{theorem}[lemma]{Theorem}
\newtheorem{corollary}[lemma]{Corollary}
\newtheorem{proposition}[lemma]{Proposition}

\theoremstyle{definition}
\newtheorem{remark}[lemma]{Remark}
\newtheorem{example}[lemma]{Example}
\newtheorem*{acknowledgement}{Acknowledgements}

\title{The heart of a convex body}

\author[Brasco]{Lorenzo Brasco} 
\address{Laboratoire d'Analyse, Topologie, Probabilit\'es, Aix-Marseille Universit\'e, CMI 39, Rue Fr\'ed\'eric Joliot Curie, 13453 Marseille Cedex 13, France} 
\email{lorenzo.brasco@univ-amu.fr}
\author[Magnanini]{Rolando Magnanini}
\address{Dipartimento di Matematica ``U. Dini'', Universit\` a di Firenze, viale Morgagni 67/A, 50134 Firenze, Italy}
\email{magnanin@math.unifi.it}

\keywords{Hot spot, eigenfunctions, convex bodies, Fraenkel asymmetry.}

\begin{document}

\maketitle

\begin{abstract}
We investigate some basic properties of the {\it heart} $\heartsuit(\mathcal{K})$ of a convex set $\mathcal{K}.$
It is a subset of $\mathcal{K},$ whose definition is based on mirror reflections of 
euclidean space, and is a non-local object. The main motivation of our interest for
$\heartsuit(\mathcal{K})$ is that this gives an estimate of the location of the hot spot in a convex heat conductor with
boundary temperature grounded at zero.
Here, we investigate on the relation between $\heartsuit(\mathcal{K})$ and the mirror symmetries of $\mathcal{K};$ 
we show that $\heartsuit(\mathcal{K})$ contains many (geometrically and phisically) relevant points of $\mathcal{K};$
we prove a simple geometrical lower estimate for the diameter of $\heartsuit(\mathcal{K});$ 
we also prove an upper estimate for the area of $\heartsuit(\mathcal{K}),$ when $\mathcal{K}$ is a triangle.
\end{abstract}

\section{Introduction}

\begin{figure}
\label{fig1}
\includegraphics[scale=.44]{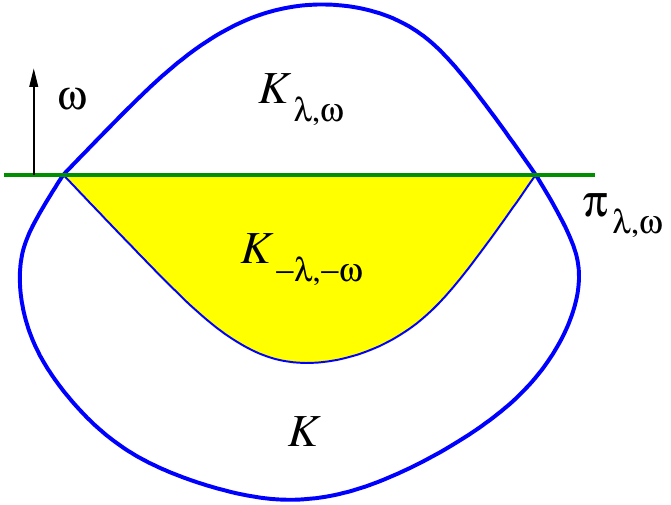}
\caption{The sets $\cK_{\la,\om}$ and $\cK_{-\la,-\om}$.}
\end{figure}

Let $\mathcal{K}$ be a convex body in the euclidean space $\mathbb{R}^N,$ that is $\mathcal{K}$ is a compact convex
set with non-empty interior.  In \cite{BMS} we defined the {\it heart} $\ccK$ of $\cK$ as follows.
Fix a unit vector $\om\in\SSN$ and a real number $\la;$ for each point $x\in\RN,$ let $T_{\la,\om}(x)$ denote the reflection of $x$ in the hyperplane $\pi_{\la,\om}$ of equation
$\langle x,\omega\rangle=\la$ (here, $\langle x,\omega\rangle$ denotes the usual scalar product of vectors in $\RN$); 
then set 
$$
\cK_{\la,\om}=\{x\in\cK: \langle x,\omega\rangle\ge\la\}
$$
(see Figure \ref{fig1}). The heart of $\cK$ is thus defined as 
$$
\ccK=\bigcap_{\om\in\SSN}\{\cK_{-\la,-\om}: T_{\la,\om}(\cK_{\la,\om})\subset\cK \}.
$$

Our interest in $\ccK$ was motivated in \cite{BMS} in connection to the 
problem of locating the (unique) point of maximal temperature --- the {\it hot spot} --- in a convex heat conductor with
boundary temperature grounded at zero. There, by means of A.D. Aleksandrov's reflection principle, we showed that $\ccK$ must contain the hot spot at each time and must also contain the maximum point of the first Dirichlet eigenfunction of the Laplacian, which is known to 
control the asymptotic behaviour of temperature for large times. By the same arguments,
we showed in \cite{BMS} that $\ccK$ must also contain the maximum point of positive solutions
of nonlinear equations belonging to a quite large class. 
By these observations, the set $\ccK$ can be viewed as a geometrical means to estimate the positions of 
these important points.

Another interesting feature of $\ccK$ is the non-local nature of its definition. 
We hope that the study of $\ccK$ can help, in a relatively simple setting, to develop techniques that may be
useful in the study of other objects and properties of non-local nature, which have lately raised 
interest in the study of partial differential equations. 
\par
A further reason of interest is that the shape of $\ccK$ seems 
to be related to the mirror symmetry of $\cK.$ By means of a numerical algorithm, developed in \cite{BMS}, 
that (approximately) constructs $\ccK$ for any given convex polyhedron $\cK,$ one can observe that relationship --- and other features of $\ccK$ --- and raise some questions. 
\begin{enumerate}
\item We know that, if $\cK$ has
a hyperplane of symmetry, then $\ccK$ is contained in that hyperplane; is the converse true?
\item How small $\ccK$ can be? Can we estimate from below the
ratio between the diameters of $\ccK$ and $\cK$? 
\item How big $\ccK$ can be? Can we estimate from above the ratio between the volumes of $\ccK$
and $\cK$?
\end{enumerate}
The purpose of this note is to collect results in that direction.
\par
In Section \ref{sec3}, we give a positive answer to question (i) (see Theorem \ref{teo:simmetry}). 

\par
In Section \ref{sec2}, we start by showing that many relevant points related to a convex set
lie in its heart. For instance, we shall prove that, besides the center of mass $M_\mathcal{K}$ of $\cK$ (as seen in \cite[Proposition 4.1]{BMS}),  
$\ccK$ must also contain the center $C_\cK$ of the smallest ball containing $\cK$ --- the so-called
{\it circumcenter} ---and
the center of mass of the set of all the {\it incenters} of $\cK,$ 
\[
\mathcal{M}(\cK)=\left\{x\in \cK\, :\, \mathrm{dist}(x,\partial \cK)=r_\cK\right\};
\]
here $r_\cK$ is the {\it inradius} of $\mathcal{K}$, i.e. the radius of the largest balls contained in $\cK$. 
This information gives a simple estimate from below of the diameter of $\ccK,$ thus partially answering to question (ii)
(see Theorem \ref{th:diameter}).
\par
By further exploring in this direction, we prove that $\ccK$ must also contain the points
minimizing each $p$-moment of the set $\cK$ (see Subsection \ref{subs:pmoment})
and other more general moments associated to $\ccK.$ As a consequence of this general
result, we relate $\ccK$ to a problem in spectral optimization considered in \cite{HKK} and 
show that $\ccK$ must contain the center of a ball realizing the 
so-called {\it Fraenkel asymmetry} (see \cite{HHW} and Section 3 for a definition).
\par
Finally, in Section \ref{sec4}, we begin an analysis of problem (iii). 
Therein, we discuss the shape optimization problem  \eqref{shape} and prove that
an optimal shape does not exist in the subclass of triangles. 
In fact, in Theorem \ref{th:triangle} we show that
$$
|\ccK|<\frac{3}{8}\,|\cK| \quad \mbox{for every triangle $\cK;$} 
$$
the constant $3/8$ is not attained but is only approached by choosing a sequence of obtuse triangles.

\section{Dimension of the heart and symmetries}
\label{sec3}

Some of the results in Subsection \ref{sub:folding}
where proved in \cite{BMS} but, for the reader's convenience, 
we reproduce them here.

\subsection{Properties of the maximal folding function}
\label{sub:folding}
The {\it maximal folding function} $\mathcal{R}_\mathcal{K}:\SSN\to \mathbb{R}$ of a convex body $\mathcal{K}\subset\mathbb{R}^N$ was defined in \cite{BMS} by
\[
\cR_\cK(\om)=\min\{\la\in\RE\ :\ T_{\la,\om}(\cK_{\la,\om})\subseteq\cK\},\quad \om\in \SSN.
\]
The heart of $\cK$ can be defined in terms of $\cR_\cK:$
\begin{equation}
\label{cuore_R}
\ccK=\left\{x\in\cK\ :\ \langle x,\omega\rangle\le\cR_\cK(\omega), \mbox{ for every } \omega\in \SSN\right\}.
\end{equation}
It is important to remark that we have 
\[
\max_{x\in\ccK} \langle x,\omega\rangle\le \mathcal{R}_\mathcal{K}(\omega),\qquad \omega\in \mathbb{S}^{N-1}
\]
and in general the two terms do not coincide: in other words, $\mathcal{R}_\mathcal{K}$ does not coincide with the support function of $\ccK$ (see \cite[Example 4.8]{BMS}).
\par
\begin{lemma}
\label{lem:lsc}
The maximal folding function $\cR_\cK:\SSN\to\mathbb{R}$ is lower semicontinuous. In particular, $\cR_\cK$ attains its minimum on $\SSN.$
\end{lemma}

\begin{proof}
Fix $\la\in\RE$ and let $\om_0\in\{ \om\in\SSN : \cR_\cK(\om)>\la\}.$ Then 
$$
|T_{\la,\om_0}(\cK_{\la,\om_0})\cap(\RN\setminus\cK)|>0.
$$
The continuity of the function
$
\om\mapsto|T_{\la,\om}(\cK_{\la,\om})\cap(\RN\setminus\cK)|
$
implies that, for every $\om$ in some neighborhood of $\om_0,$ we have
that $\cR_\cK(\om)>\la.$ 
\end{proof}
\begin{remark}
In general $\mathcal{R}_\cK$ is not continuous on $\SSN$: it sufficient to take a rectangle
$\cK=[-a,a]\times[-b,b]$ and observe that we have $\cR_\cK(1,0)=0,$
since $\cK$ is symmetric with respect to the $y$-axis but, for a sufficiently small $\tht,$
$
\cR_\cK(\cos\tht,\sin\tht)=a\,\cos\tht- b\,\sin\tht,
$
so that
\[
\lim_{\tht \to 0} \mathcal{R}_\cK (\tht)=a>\mathcal{R}_\cK(1,0).
\]

Observe that here the lack of continuity of the maximal folding function is not due to the lack of smoothness of the boundary of $\cK$, but rather to the presence of non-strictly convex subsets of $\partial \cK$. 
\end{remark}

\begin{proposition}
\label{lm:baricentro}
Let $\cK\subset\mathbb{R}^N$ be a convex body and define its {\it center of mass} by
\[
M_\cK=\frac1{|\cK|}\int_{\cK} y\, dy.
\]
Then we have that
\begin{equation}
\label{positiva}
\cR_{\cK}(\om)\ge \langle M_\cK,\om\rangle,\quad \mbox{ for every } \om\in \SSN,
\end{equation}
and the equality sign can hold for some $\om\in\SSN$ if and only if $\cK$ is $\om-$symmetric, i.e. if 
$T_{\la,\om}(\cK)=\cK$ for some $\la\in\RE.$
In particular, $M_\cK\in\ccK.$
\end{proposition}
\begin{proof}
Set $\la=\cR_\cK(\om)$ and consider the set $\Om=\cK_{\la,\om} \cup T_{\la,\om}(\cK_{\la,\om});$ 
by the definition of center of mass and since $\Om$ is symmetric with respect to the hyperplane $\pi_{\la,\om},$  
we easily get that
\[
\begin{split}
|\cK|\, \left[\cR_\cK(\omega)-\langle M_\cK,\om\rangle\right]&=\int_\cK \big[\cR_\cK(\om)-\langle y,\om\rangle\big]\, dy\\
&=\int_{\cK\setminus\Om} \big[\cR_\cK(\om)-\langle y,\om\rangle\big]\, dy.
\end{split}
\]
Observe that the last integral contains a positive quantity to be integrated on the region $\cK\setminus\Om$: this already shows \eqref{positiva}. Moreover, the same identity implies:
\[
\cR_\cK(\om)-\langle M_\cK,\om\rangle=0\quad \Longleftrightarrow \quad |\cK\setminus\Om|=0,
\]
and the latter condition is equivalent to say that $\cK$ is $\omega-$symmetric. Finally, 
by combining \eqref{positiva} with the definition \eqref{cuore_R} of $\ccK,$ it easily follows $M_\cK\in\ccK.$
\end{proof}

\subsection{On the mirror symmetries of $\cK$}

A first application of Proposition \ref{lm:baricentro} concerns 
the relation between $\ccK$ and the mirror symmetries of $\cK.$ 
\begin{theorem}
\label{teo:simmetry}
Let $\cK\subset\RN$ be a convex body. 
\begin{enumerate}
\item[(i)]
If there exist $k$ $(1\le k\le N)$ independent directions $\om_1,\dots,\om_k\in\SSN$ such that
$\cK$ is $\om_j$-symmetric for $j=1,\dots, k,$ then $\cR_\cK(\om_j)=\langle M_\cK,\omega_j\rangle$ for $j=1,\dots, k$
and
$$
\ccK\subseteq\bigcap_{j=1}^k \pi_{\cR_\mathcal{K}(\om_j),\om_j}.
$$
In particular, the co-dimension of $\ccK$ is at least $k.$
\item[(ii)]  If $\ccK$ has dimension $k$ $(1\le k\le N-1),$ then there exists at least a direction
$\theta\in\SSN$ such that
$\cR_\cK(\theta)=\langle M_\cK,\theta\rangle$, and $\cK$ is $\theta$-symmetric.
\end{enumerate}
\end{theorem}
\begin{proof}
(i) The assertion is a straightforward consequence of Proposition \ref{lm:baricentro}.
\vskip.2cm
(ii) By Lemma \ref{lem:lsc}, the function $\om\mapsto\cR_\cK(\om)-\langle M_\cK,\om\rangle$ attains
its minimum for some $\theta\in\SSN.$ 
Set $r=\cR_\cK(\theta)-\langle M_\cK,\theta\rangle$ and suppose that $r>0.$
\par
Then, for every $x\in B(M_\cK,r),$ we have that
$$
\langle x,\om\rangle=\langle M_\cK,\om\rangle+\langle x-M_\cK,\om\rangle<\langle M_\cK,\om\rangle+r\le
\cR_\cK(\om)-r+r=\cR_\cK(\om),
$$ 
for every $\om\in\SSN,$ and hence $x\in\ccK$ by \eqref{cuore_R}. Thus, $B(x,r)\subset\ccK$ 
--- a contradiction to the fact that $1\le k\le N-1.$ Hence, 
$\cR_\cK(\theta)=\langle M_\cK,\theta\rangle$ and $\cK$ is $\theta$-symmetric by Proposition \ref{lm:baricentro}.
\end{proof}

\begin{remark}
It is clear that the dimension of the heart only gives information on the {\it minimal} number of symmetries of a convex body: the example of a ball is quite explicative.
\par\noindent
We were not able to prove the following result:
\vskip.2cm
{\it if $\ccK$ has co-dimension $m$ $(1\le m\le N),$ then
there exist at least $m$ independent directions $\theta_1,\dots, \theta_N\in\SSN$ such that 
$\cR_\cK(\theta_j)=\langle M_\cK, \theta_j\rangle,$ $j=1,\dots, m,$ and $\cK$ is
$\theta_j$-symmetric for $j=1,\dots, m.$}
\vskip.2cm\noindent
We leave it as a conjecture.
\end{remark}

\section{Relevant points contained in the heart}
\label{sec2}

In this section, we will show that many relevant points of a convex set are contained in its heart 
(e.g. the incenter and the circumcenter, besides the center mass). This fact will give us a means to estimate from below
the diameter of the heart.

\subsection{An estimate of the heart's diameter}

\begin{proposition}
\label{prop:circo}
Let $\cK\subset\RN$ be a convex body.
Then its circumcenter $C_\cK$ belongs to $\ccK.$
\end{proposition}
\begin{proof}
Suppose that $C_\cK\notin\ccK$ and let $q\in\ccK$ be its (unique) projection on the set $\ccK$; then, define:
$$
H^+=\{x\in\RN: \langle x-q,C_\cK-q\rangle>0\} \qquad \mbox{ and } \qquad \Pi=\pa H^+.
$$
Since $C_\cK\in\cK,$ then $H^+\cap\cK$ is non-empty and its reflection in $\Pi$ is contained in $\cK,$ by the definition of $\ccK,$ and hence is contained in $B_R(C_\cK).$
Thus, $\cK$ must be contained in the set $B_R(C_\cK)\cap B_R(2q-C_\cK)$.
\par
This is a contradiction, since the smallest ball containig $\cK$ would have
a radius that is at most $\sqrt{R^2-|C_\cK-q|^2}<R$.
\end{proof}
We now consider the {\it incenters} of $\cK$: these are the centers of the balls of largest radius $r_\cK$ inscribed in $\cK$. Needless to say, a convex body may have many incenters. We start with the simpler case of a convex body with a unique incenter.
\begin{proposition}
\label{prop:in}
Let $\cK\subset\RN$ be a convex body; if its incenter $I_\cK$ is unique, then 
$I_\cK\in\ccK.$ In particular, $I_\cK\in\ccK$ if $\cK$ is strictly convex.
\end{proposition}
\begin{proof}
Consider the unique maximal ball $B(I_\cK, r_\cK)$ inscribed in $\cK$  and suppose that 
$I_\cK\notin\ccK;$  this implies that there exists $\omega\in\SSN$ such that 
\[
\mathcal{R}_\cK(\omega)-\langle I_\cK,\omega\rangle<0.
\]
Set $\la=\cR_\cK(\om)$ and define $I'_\cK=T_{\la,\omega}(I_\cK);$ then $I'_\cK\not=I_\cK$ and $\langle I'_\cK,\omega\rangle<\langle I_\cK,\omega\rangle.$
Now, the half-ball
$B^+=\{x\in B(I_\cK, r_\cK)\, :\, \langle x,\omega\rangle\ge \la)\}$
and its reflection $T_{\la,\om}(B^+)$ in the hyperplane $\pi_{\la,\om}$ 
are contained in $\cK,$ since $B^+$ is contained in the maximal cap $\cK_{\la,\om}.$
\par
This fact implies in particular that the reflection of the whole ball $B$ is contained in $\cK$: but the latter is still a maximal ball of radius $r_\cK$, with center $I'_\cK$ different from $I_\cK$.  This is a contradiction, 
since it violates the assumed uniqueness of the incenter.
\end{proof}

To treat the general case, we need the following simple result.

\begin{lemma}
Let $\cK\subset\RN$ be a convex body and let us set
\[
\cM(\cK)=\left\{x\in \cK\, :\, \mathrm{dist}(x,\partial \cK)=r_\cK\right\}.
\]
Then $\cM(\cK)$ is a closed convex set with $|\cM(\cK)|=0;$ in particular, the dimension of $\cM(\cK)$ is at most $N-1$.
\end{lemma}
\begin{proof}
The quasi-convexity\footnote{This means that the superlevel sets of the function are convex.} of $\dist(x,\partial\cK)$
(being $\cK$ convex) immediately implies that $\cM(\cK)$ is convex.   
For the reader's convenience, here we give a proof anyway.
Let us take $x,z\in\cM(\cK)$ two distinct points, then by definition of inradius we have
\[
B(x,r_\cK)\cup B(z,r_\cK)\subset \cK.
\]
Since $\cK$ is convex, it must contain the convex hull of $B(x,r_\cK)\cup B(z,r_\cK),$ as well: hence, for every $t\in[0,1]$ we have
\[
B((1-t)\,x+t\,z,r_\cK)\subset\cK,
\]
that is $(1-t)\,x+t\,w\in\cM(\cK)$, which proves the convexity of $\cM(\cK)$.
\par
Now, suppose that $|\cM(\cK)|>0;$ since $\cM(\cK)$ is convex, it must contain a ball $B(x,\varrho)$. 
The balls of radius $r_\cK$ having centers on $\pa B(x,\varrho/2)$ are all contained in $\mathcal{K}$, so that their whole union is contained in $\cK$ as well. Observing that this union is given by $B(x,\varrho/2+r_\mathcal{K})$, we obtain the desired contradiction, since we violated the maximality of $r_\mathcal{K}$.
\end{proof}

\begin{proposition}
\label{prop:ins}
Let $\cK\subset\mathbb{R}^N$ be a convex body and let us suppose that $\cM(\cK)$ has dimension $k\in \{1,\dots,N-1\}$. 
\par
Then the center of mass of $\cM(\cK),$ defined by
\[
I_\cM=\frac{\displaystyle\int_{\cM(\cK)} y\, d\cH^k(y)}{\cH^k(\cM(\cK))},
\]
belongs to $\ccK$. Here, $\cH^k$  denotes the standard $k-$dimensional Hausdorff measure.
\end{proposition}
\begin{proof}
The proof is based on the observation that 
\begin{equation}
\label{complanari}
\mathcal{R}_{\cM(\cK)}(\omega)\le \mathcal{R}_\cK(\omega),\qquad \omega\in \SSN,
\end{equation}
where $\mathcal{R}_{\cM(\cK)}$ is the maximal folding function of $\cM(\cK)$, thought as a subset of $\mathbb{R}^N$. Assuming \eqref{complanari} to be true, we can use the definition of heart and Proposition \ref{lm:baricentro} to obtain
that the center of mass $I_\cM$ belongs to $\heartsuit(\cM(\cK))$ and hence to $\ccK,$
which would conclude the proof.
\par
Now, suppose by contradiction that there is an $\om\in\SSN$ such that
$\cR_{\cM(\cK)}(\om)> \cR_\cK(\om),$
and set $\la=\cR_\cK(\om),$ as usual. Then, there exists $x\in\cM(\cK)$ with $\langle x,\omega\rangle\ge\la$ such that its reflection $x^\la=T_{\la,\om} (x)$ in the hyperplane $\pi_{\la,\om}$ falls outside $\cM(\cK)$: this would imply in particular that $B(x^\lambda,r_\mathcal{K})\not\subset\mathcal{K}$. Observe that by definition of $\cM(\cK)$, the ball $B(x,r_\mathcal{K})$ lies inside $\cK$, so that the cap
$B(x,r_\mathcal{K})\cap\{\langle y,\om\rangle\ge \la\}$ 
is reflected in $\cK;$ thus, as before, we obtain that the union of this cap and its reflection is contained in $\cK$. This shows that the ball $B(x^\la,r_\mathcal{K})$ is contained in $\cK$, thus giving a contradiction. 
\end{proof}

The result here below follows at once from Propositions \ref{lm:baricentro}, \ref{prop:circo} and \ref{prop:ins}.

\begin{theorem}
\label{th:diameter}
Let $\cK\subset\RN$ be a convex body, then
$$
\diam[\ccK]\ge \max(|M_\cK-C_\cK|,|C_\cK-I_\cM|,|I_\cM-M_\cK|).
$$
\end{theorem}

\begin{remark}
Notice that, when $\ccK$ degenerates to a single point, then clearly
\[
\ccK=\{M_\cK\}=\{I_\cK\}=\{C_\cK\}.
\]
\par
Needless to say, it may happen that the three points $M_\cK,I_\cK$ and $C_\cK$ coincide, but $\ccK$ is not a point. 
For example, take an ellipse parametrized in polar coordinates as
\[
E=\left\{(\varrho,\vartheta)\, :\, 0\le \varrho\le \sqrt{a^2\cos^2\vartheta+b^2\, \sin^2\vartheta}, \
\tht\in[-\pi,\pi]\right\}
\]
and a $\pi$-periodic, smooth non-negative function $\eta$ on $[-\pi,\pi],$
having its support in two small neighborhoods of $-3\pi/4$ and $\pi/4.$ 
Then, 
if $\ep$ is sufficiently small, the deformed set
\[
E_\ep=\left\{(\varrho,\tht)\, :\, 0\le \varrho\le \sqrt{a^2\cos^2\tht+b^2\, \sin^2\tht}-\ep\,\eta(\tht), \
\tht\in[-\pi,\pi]\right\}
\]
is still convex and centrally symmetric. Moreover, it is easy
to convince oneself that $\{M_{E_\ep}\}=\{I_{E_\ep}\}=\{C_{E_\ep}\}=\{(0,0)\}$, whereas, by Theorem \ref{teo:simmetry}, $\heartsuit(E_\ep)$ is not a point, since $E_\ep$ has no mirror symmetries.
\end{remark}

\subsection{The p-moments of $\cK$ and more}
\label{subs:pmoment}
We recall that the point $M_\cK$ can also be characterized as the unique point in $\cK$ which minimizes the function
\[
x\mapsto\int_\cK |x-y|^2\, dy,\qquad x\in\mathcal{K},
\]
that can be viewed as the moment of inertia (or $2$-moment) of $\cK$ about the point $x.$
In this subsection, we will extend the results of Subsection 3.1 to 
more general moments of $\cK,$ that include as special cases the $p$-moments $\int_\cK |x-y|^p\, dy.$ 
\par
We first establish a preliminary lemma.

\begin{lemma}
\label{lm:conticino}
Let $\fhi, \psi:[0,\infty)\to\mathbb{R}$ be, respectively, an increasing and a decreasing function
and suppose that $\psi$ is also integrable in $[0,\infty).$ Define the two functions 
\[
F(t)=\int_a^b \varphi(|s-t|)\, ds,\quad t\in \RE,
\]
and
\[
G(t)=\int_{-\infty}^a \psi(t-s)\, ds+\int_b^\infty \psi(s-t)\, ds,\quad t\in \RE,
\]
where $a$ and $b$ are two numbers with $0\le a\le b.$
\par
Then, both $F$ and $G$ attain their minimum at the 
midpoint $(a+b)/2$ of $[a,b].$ If $\fhi$ is strictly increasing (resp. $\psi$ is strictly decreasing),
then the minimum point of $F$ (resp. $G$)  is unique. 
\end{lemma}
\begin{proof}
(i) Let us introduce a primitive of $\varphi$, i.e.
\[
\Phi(t)=\int_{0}^t \varphi(s)\, ds;
\]
observe that this is a convex function, since its derivative is increasing. Now, we can write
\[
F(t)=\int_a^t \varphi(t-s)\, ds+\int_t^b \varphi(s-t)\, ds=\Phi(t-a)+\Phi(b-t),
\]
for $t\in [a,b]$, 
and notice that both functions on the right-hand side are convex, thus $F$ is convex in $[a,b]$, as well. 
It is not difficult to see that the graph of $F$ in $[a,b]$ is symmetric with respect to the line $t=(a+b)/2$, indeed for every $t\in [a,b]$ we easily get that
\[
F(a+b-t)=F(t).
\] 
This shows that $F$ attains its minimum in $[a,b]$ at the midpoint. Clearly, if $\fhi$ is strictly increasing, then $F$ is strictly convex and the minimum point is unique.
\par 
Finally, we observe that $F$ attains at $(a+b)/2$ also the global minimum on $\RE$, since
$F$ is clearly decreasing in $(-\infty,a]$ and increasing on $[b,+\infty)$.
\par
(ii) It is enough to rewrite the function $G$ as follows
\[
G(t)=\int_{\mathbb{R}} \psi(|t-s|)\, ds-\int_{a}^b \psi(|t-s|)\, ds,\qquad t\in[a,b];
\]
then we may notice that the first term is a constant, while the second one behaves as $F$, thanks to the first part of this lemma, since the function $-\psi$ is increasing.
\end{proof}
\par 
The following result generalizes one in \cite{OH} in the case of $\ccK$ convex.
\begin{theorem}
\label{teo:pmoment}
Let $\cK\subset\mathbb{R}^N$ be a convex body and let $\fhi:[0,\infty)\to\RE$ be an increasing function.
Define the function 
$$
\mu_\fhi(x)=\int_\cK \fhi(|x-y|)\,dy,\quad x\in\RN,
$$
and the set 
\[
\mathbf{m}(\mu_\fhi)=\{x\in\mathbb{R}^N\, :\, \mu_\varphi(x)=\min \mu_\varphi\}.
\]
Then
\begin{equation}
\label{intersezione}
\mathbf{m}(\mu_\varphi)\cap\ccK\not=\varnothing.
\end{equation}
\end{theorem}
\begin{proof}
We shall refer to $\mu_\fhi(x)$ as the {\it $\fhi$-moment} of $\cK$ about the point $x.$ 
First of all, we observe that $\mu_\varphi$ is lower semicontinuous thanks to Fatou's Lemma and that 
\[
\inf_{\mathcal{K}} \mu_\varphi=\inf_{\mathbb{R}^N} \mu_\varphi,
\]
so that the minimum of $\mu_\varphi$ is attained at some point belonging to $\mathcal{K}$, i.e. $\emptyset\not=\mathbf{m}(\mu_\varphi)\subset\mathcal{K}$.
\vskip.2cm\noindent
We first prove \eqref{intersezione} when $\fhi$ is strictly increasing. Let $x\in\mathcal{K}$ be a minimum point of $\mu_\fhi$.
If $x\notin\ccK$, then there exists a direction $\omega\in\mathbb{S}^{N-1}$ such that $\mathcal{R}_\mathcal{K}(\omega)<\langle x,\omega\rangle$. We set for simplicity $\lambda=\mathcal{R}_\mathcal{K}(\omega)$ and consider the hyperplane $\pi=\pi_{\la,\om}$, so that
\[
x\in \mathcal{K}_{\lambda,\omega} \qquad \mbox{ and } \qquad\mathcal{T}_{\lambda,\omega}(\mathcal{K}_{\lambda,\omega})\subset\mathcal{K}.
\] 
Modulo a rotation, we can always assume that $\omega=(1,0,\dots,0)$ and the hyperplane $\pi$ has
the form $\{x\in\mathbb{R}^N\, :\, x_1=\la\}$. 
\par
Now, we define the symmetric set $\Om=\cK_{\la,\om}\cup\cT_{\la,\om}(\cK_{\la,\om})$, which 
can be witten as
$$
\Om=\{(y',y_1)\in\cK: y'\in\cK\cap\pi, \la-a(y')\le y_1\le \la+a(y')\}.
$$
Consider the projection $z$ of $x$ in the hyperplane $\pi$ and observe that $z\in \Omega$, thanks to the convexity of $\mathcal{K}$. 
For $y\in \cK\setminus\Om$,
we have that $|z-y|<|x-y|$. Thus
\begin{equation}
\label{integ}
\int_{\cK\setminus\Om} \fhi(|x-y|)\,dy\ge\int_{\cK\setminus\Om} \fhi(|z-y|)\,dy,
\end{equation}
since $\fhi$ is increasing.  Moreover, by Fubini's theorem, we compute:
\[
\int_{\Om} \fhi(|x-y|)\,dy=
\int_{\cK\cap\pi}\Biggl\{\int_{\la-a(y')}^{\la+a(y')}
\fhi\left(\sqrt{|x'-y'|^2+|x_1-y_1|^2}\right)\,dy_1\Biggr\}\,dy'.
\]
We now apply Lemma \ref{lm:conticino} with the choice
\(t\mapsto\fhi\left(\sqrt{|x'-y'|^2+t^2}\right)\), which is a strictly increasing function. Thus, we can infer that the last integral 
is strictly larger than  
\[
\int_{\cK\cap\pi}\Biggl\{\int_{\la-a(y')}^{\la+a(y')}
\fhi\left(\sqrt{|x'-y'|^2+|\la-y_1|^2}\right)\,dy_1\Biggr\}\,dy'=
\int_{\Om} \fhi(|z-y|)\,dy.
\]
With the aid of \eqref{integ}, we then conclude that $\mu_\fhi(x)>\mu_\fhi(z)$, which gives
a contradiction. We observe in passing that we have proved something more, namely,
we showed that $\mathbf{m}(\mu_\fhi)\subseteq\ccK$.
\vskip.2cm\noindent
If $\fhi$ is only increasing, we approximate it by the following sequence of strictly increasing
functions 
\[
\fhi_n(t)=\varphi(t)+\frac{1}{n}\, t^2,\qquad t\in[0,\infty),\ n\in\mathbb{N}.
\]
Let $x_n\in\mathbf{m}(\mu_{\fhi_n})\subseteq\ccK$ be a sequence of minimizers of $\mu_{\varphi_n}$, then by compactness of $\heartsuit(\mathcal{K})$ they convergence (up to a subsequence) to a point $x_0\in\heartsuit(\mathcal{K})$. For every $x\in\mathbb{R}^N$, by Fatou's Lemma we have
\[
\mu_\varphi(x)=\lim_{n\to\infty} \mu_{\varphi_n}(x)\ge \liminf_{n\to\infty} \mu_{\varphi_n}(x_n)\ge \mu_\varphi(x_0),
\]
which implies that $x_0\in\mathbf{m}(\mu_\varphi)$. This concludes the proof.
\end{proof}
The analogous of Theorem \ref{th:diameter} is readily proved.
\begin{theorem}
\label{teo:diamgen}
Let $\cK$ be a convex body. Then the convex hull of the set
$$
\bigcup\left\{\mathbf{m}(\mu_\fhi)\cap\ccK: \fhi \mbox{ is increasing on $[0,\infty)$}\right\},
$$ 
is contained in $\ccK.$
\end{theorem}
\begin{remark}
By similar arguments, we can prove that, if $\psi$ is decreasing and the function
$$
\nu_\psi(x)=\int_{\RN\setminus\cK} \psi(|x-y|)\,dy,
$$
is finite for every $x\in\cK$, then $\mathbf{m}(\nu_\psi)\cap\ccK\not=\varnothing.$
\end{remark}
\par
Particularly interesting are the cases where $\fhi(t)=t^p$ with $p>0$ 
and $\psi(t)=t^{-p}$ with $p>N.$ 
\begin{corollary} 
\label{coro:p}
Consider the functions 
$$
\mu_p(x)=\int_\cK |x-y|^p\,dy,\quad x\in\cK,
$$
for $p>0$ or
$$
\nu_p(x)=\int_{\RN\setminus\cK} |x-y|^{-p}\,dy,\quad x\in\cK,
$$
for $p>N.$ 
\par
Then, their minimum points belong to $\ccK.$
\end{corollary}
\begin{remark}
Propositions \ref{lm:baricentro}, \ref{prop:circo} and \ref{prop:in} can be
re-proved by means of Corollary \ref{coro:p} by choosing $p=2$ or, respectively,
by taking limits as $p\to\infty.$ 
\par
Notice, in fact, that
$$
\lim_{p\to\infty} \mu_p(x)^{1/p}=\max_{y\in\cK}|x-y|
$$
and
$$
\lim_{p\to\infty} \nu_p(x)^{-1/p}=\min_{y\in\cK}|x-y|.
$$
Hence the {\it circumradius} $\rho_\cK$ and inradius $r_\cK$ are readily obtained as 
$$
\rho_\cK=\min_{x\in\cK}\,\max_{y\in\pa\cK} |x-y|=\lim_{p\to+\infty} \min_{x\in\cK}\,\mu_p(x)^{1/p},
$$
and 
$$ 
r_\cK=\max_{x\in\cK}\, \min_{y\in\pa\cK} |x-y|=\lim_{p\to+\infty} \max_{x\in\cK}\,\nu_p(x)^{-1/p}.
$$
These observations quite straightforwardly imply that $C_\cK$ and $I_\cK$ belong to $\ccK.$
\par
A final remark concerns the case $p=0.$ It is well-known that
$$
\lim_{p\to 0^+} \left(\frac{\mu_p(x)}{|\cK|}\right)^{1/p}=\exp\left\{\int_\cK \log|x-y|\,\frac{dy}{|\cK|}\right\}=
\exp\{\mu_{\log}(x)/|\cK|\},
$$
that can be interpreted as the {\it geometric mean} of the function $y\mapsto |x-y|$ on $\cK;$
needless to say, the set of its minimum points intersects $\ccK.$
\end{remark}

\subsection{On Fraenkel's asymmetry}

\begin{lemma}
\label{lm:ga}
Let $\cK\subset\RN$ be a convex body. For $r>0,$ define the function
\[
\ga(x)=|\cK\cap B(x,r)|,\qquad x\in\RN.
\]
Then $\ga$ is log-concave and, if $\mathbf{M}(\ga)=\{x\in\RN, :\, \ga(x)=\max \ga\},$ then
\begin{equation}
\label{massimoG}
\mathbf{M}(\ga)\cap \heartsuit(\mathcal{K})\not=\varnothing.
\end{equation}
\end{lemma}
\begin{proof}
The log-concavity of $G$ is a consequence of Pr\'ekopa--Leindler's inequality, that
we recall here for the reader's convenience:
let $0<t<1$ and let $f,g$ and $h$ be nonnegative integrable functions on $\RN$ satisfying
\begin{equation}
\label{hpl}
h((1-t)\, x+t\, y)\ge f(x)^{1-t}\, g(y)^t,\qquad \mbox{ for every }x,y\in\RN;
\end{equation}
then
\begin{equation}
\label{pl}
\int_{\mathbb{R}^N} h(x)\, dx\ge \left(\int_{\mathbb{R}^N} f(x)\,dx\right)^{1-t}\left(\int_{\mathbb{R}^N}g(x)\, dx\right)^t.
\end{equation}
(For a proof and a discussion on the links between \eqref{pl} and the Brunn-Minkowski inequality, the reader is referred to \cite{Ga}.)
\par
Indeed, we pick two points $z ,w\in\mathbb{R}^N$ and a number $t\in(0,1),$  
and apply Pr\'ekopa--Leindler's inequality to the triple of functions
\[
f=1_{\cK\cap B(z,r)},\quad g=1_{\cK\cap B(w,r)},\quad h=1_{\cK\cap B((1-t)\,z+t\, w,r)};
\]
then, \eqref{hpl} is readily satisfied.
Thus, \eqref{pl} easily implies
\[
\begin{split}
\ga((1-t)\, z+t\, w)=&|\cK\cap B((1-t)\,z+t\,w,r)|\ge \\
&|\cK\cap B(z,r)|^{1-t}\, |\cK\cap B(w,r)|^t=
\ga(z)^{1-t}\,\ga(w)^t,
\end{split}
\]
and, by taking the logarithm on both sides, we get
the desired convexity. 
A straightforward consequence is that the set $\mathbf{M}(\ga)$
is convex.
\par
Once again, the validity of \eqref{massimoG} will be a consequence of the inequality
\begin{equation}
\label{G}
\cR_{\mathbf{M}(\ga)}(\om)\le\cR_\cK(\om),\quad \mbox{ for every }\om\in\SSN.
\end{equation}
\par
By contradiction: let us suppose that there exist $\om\in\SSN$ and $x\in\mathbf{M}(\ga)$ such that
\[
\cR_\cK(\om)<\cR_{\mathbf{M}(\ga)}(\om)\le \langle x,\om\rangle.
\] 
In particular, this implies that the point $x^\la=T_{\la,\om}(x),$ with $\la=\cR_\cK(\om),$
does not belong to $\mathbf{M}(\ga)$ -- i.e. the reflection of $x$ with respect to the hyperplane $\pi_{\la,\om}$ falls outside $\mathbf{M}(\ga)$.
\par
We set for brevity $B=B(x,r)$ and $B^\la=B(x^\la,r),$ 
and we again consider the $\om-$symmetric set
$\Om=\cK_{\la,\om}\cup T_{\la,\om}(\cK_{\la,\om})\subseteq\cK.$
Then, observe that
\[
\begin{split}
B\cap (\cK\setminus \Om)&=\{x\in B\, :\, \langle x,\om\rangle<\la\}\cap (\cK\setminus\Om)\\
&\subseteq (B\cap B^\la)\cap (\cK\setminus\Om)\subseteq B^\la\cap (\cK\setminus\Om),
\end{split}
\]
which implies that $|B\cap (\cK\setminus \Om)|\le |B^\la\cap(\cK\setminus \Om)|.$
Also, notice that since $\Om$ is symmetric in the hyperplane $\pi_{\la,\om}$ and $B^\la=T_{\la,\om}(B)$, we have
that $|B^\la\cap\Om|=|B\cap\Om|.$
\par
By using these informations and the maximality of $x$, we can infer that
\[
\begin{split}
\ga(x^\la)=|\cK\cap B^\la|&=|\Om\cap B^\la|+|(\cK\setminus\Om)\cap B^\la|\\
&\ge |\Om\cap B|+|(\cK\setminus\Om)\cap B|=|\Om\cap B|=\ga(x),
\end{split}
\]
that is $x^\la$ is also a maximum point, i.e. $x^\la\in\mathbf{M}(\ga)$
--- a contradiction.
\end{proof}
As a consequence of this lemma, we obtain a result concerning the so-called 
{\it Fraenkel asymmetry} of $\cK:$ 
\[
\cA(\cK)=\min_{x\in \RN}\frac{|\cK\De B(x,r^*_\cK)|}{|\cK|},
\]
where $\De$ denotes the symmetric difference of the two sets and 
the radius $r^*_\cK$ is determined by $|B(x,r^*_\cK)|=|\cK|.$
This is a measure of how a set is far
from being spherically symmetric and was introduced in \cite{HHW}; we refer the reader to \cite{FMP} 
for a good account on $\cA(\cK).$
\begin{corollary}
Let $\cK\subset \RN$ be a convex body. 
Then $\cA(\cK)$ is attained for at least one ball centered at a point belonging to $\ccK$. 
\end{corollary}
\begin{proof}
It is sufficient to observe that  
$|\cK\De B(x,r^*_\cK)|=2(|\cK|-|\cK\cap B(x,r^*_\cK)|)$, since $|B(x,r^*_\cK)|=|\cK|,$ and hence
\[
\frac{|\cK\De B(x,r^*_\cK)|}{|\cK|}=2\left(1-\frac{\ga(x)}{|\cK|}\right).
\]
Thus, $\cA(\cK)$ is attained by points that maximize $\ga;$ hence, Lemma \ref{lm:ga}
provides the desired conclusion. 
\end{proof}
\begin{remark}
Observe in particular that if $\cK$ has $N$ hyperplanes of symmetry, then an optimal ball can be placed at their intersection. However, in general, even under this stronger assumption, such optimal ball is not unique. For example, 
take the rectangle $Q_\ep=[-\pi/4\ep,\pi/4\ep]\times[-\ep,\ep]$ with $0<\ep<\pi/4;$
any unit ball centered at a point in the segment $(-\pi/4\ep+1, \pi/4\ep-1)\times\{0\}$
realizes the Fraenkel asymmetry $\cA(Q_\ep)$. Thus, in general {\it it is not true that all optimal balls are centered in the heart}.
\end{remark}

\begin{remark}
The following problem in spectral optimization was considered in \cite{HKK}: given a (convex) set $\cK\subset\mathbb{R}^N$ 
and a radius $0<r<r_\cK$, find the ball $B(x_0,r)\subset\cK$ which maximizes the quantity
\[
\lambda_1(\cK\setminus B(x,r))
\]
as a function of $x$: here, $\lambda_1(\Omega)$ stands for the first Dirichlet-Laplacian eigenvalue of a set $\Omega$. By considerations similar to the ones used in this section and
remarks contained in \cite[Theorem 2.1]{HKK}, it can be proved that $x_0\in\ccK.$
\end{remark}

\section{Estimating the volume of the heart}
\label{sec4}

In this section, we begin an analysis of the following problem in
shape optimization:
\begin{equation}
\label{shape}
\mbox{maximize the ratio } \ \frac{|\ccK|}{|\cK|} \ \mbox{ among all convex bodies $\cK\subset\RN$};
\end{equation}
solving \eqref{shape} would give an answer to question (iii) in the introduction. Since this ratio is scaling invariant, \eqref{shape} is equivalent to the following problem:
\begin{equation}
\label{diameter}
\mbox{maximize the ratio } \ \frac{|\ccK|}{|\cK|} \ \mbox{ among all convex bodies $\cK\subset [0,1]^N$};
\end{equation}
here, $[0,1]^N$ is the unit cube in $\RN.$ 
\par
We notice that the class of the competing sets in problem \eqref{diameter} 
is relatively compact in the topology induced by the {\it Hausdorff distance} (see \cite[Chapter 2]{HP})
 --- the most natural topology when one deals with the constraint of convexity. 
This fact implies, in particular, that any maximizing sequence $\{\cK_n\}_{n\in\mathbb{N}}\subset[0,1]^N$ of convex bodies converges 
-- up to a subsequence -- to a compact convex set $\mathcal{K}\subset[0,1]^N$. 
\par
However, there are two main obstructions to the existence of a maximizing set for \eqref{diameter}:
(a) in general, the limit set $\mathcal{K}$ may not be a convex body, i.e. $\mathcal{K}$ could have empty interior; 
in other words, maximizing sequences could ``collapse'' to a lower dimensional object; 
(b) it is not clear whether the shape functional $\cK\mapsto|\ccK)|$
is upper semicontinuous or not in the aforementioned topology. 
\par
The next example assures that the foreseen semicontinuity property {\it fails to be true} in general.
\begin{example}
Let $Q=[-2,2]\times[-1,1]$ and take the points 
\[
p^1_\ep=(1,1+\ep)\qquad \mbox{ and }\qquad p^2_\ep=(-2-\varepsilon,1/2),
\] 
and define
$Q_\ep$ as the convex hull of $Q\cup\{p^1_\ep,\,p^2_\ep\}.$ As $\ep$ vanishes,
$\heartsuit(Q_\ep)$ shrinks to the quadrangle having vertices 
$$
(0,0)\qquad (1,0) \qquad (1/2,1/2)\qquad \mbox{ and } \qquad (0,1/2),
$$
while clearly $\heartsuit(Q)=\{0\}$: indeed, observe that due to the presence of the new corners $p^1_\varepsilon$ and $p^2_\varepsilon$, it is no more possible to use $\{x=0\}$ and $\{y=0\}$ as maximal axis of reflection in the directions $\mathbf{e}_1=(1,0)$ and $\mathbf{e}_2=(0,1)$. In particular, we get that
\[
0=|\heartsuit(Q)|<\lim_{\varepsilon\to 0^+} |\heartsuit(Q_\varepsilon)|.
\]
\begin{figure}[h]
\label{fig2}
\includegraphics[scale=.25]{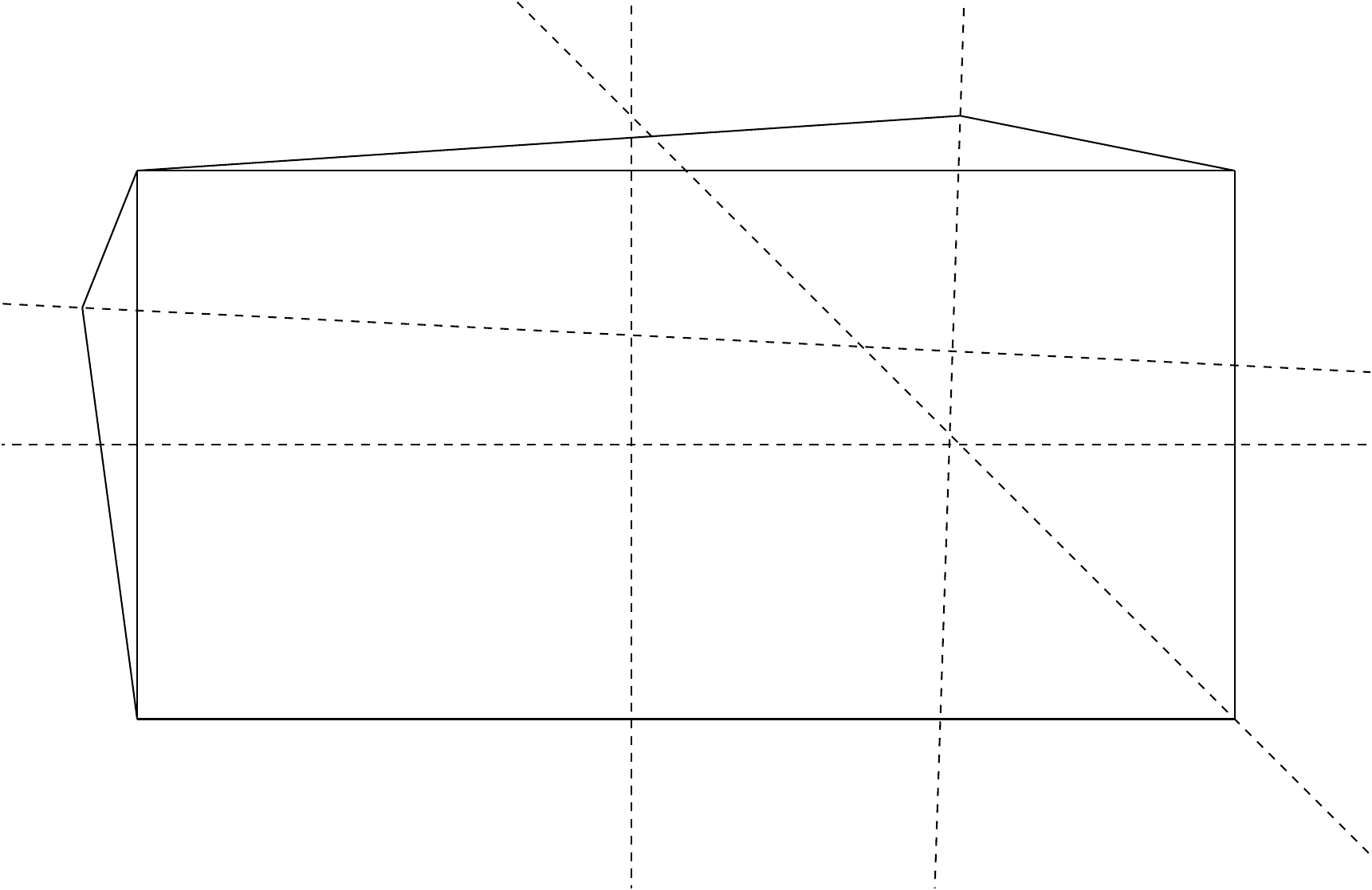}
\caption{The heart of $Q_\ep$.}
\end{figure}
\end{example}
These considerations show that the existence of a solution of \eqref{shape} is not a trivial issue.
Indeed, we are able to show that {\it an optimal shape does not exist} in the class of triangles. 
This is the content of the main result of this section.
\begin{theorem}
\label{th:triangle}
It holds that
$$
\sup\left\{ \frac{|\ccK|}{|\cK|}: \cK \mbox{ is a triangle} \right\}=\frac38,
$$
and the supremum is attained by a sequence of obtuse triangles.
\end{theorem}

The proof of Theorem \ref{th:triangle} is based on the following lemma, in which we 
exactly determine $\ccK,$ when $\cK$ is a triangle. 

\begin{lemma}
\label{th:hearttriangle}
Let $\cK$ be a triangle.  Then the following assertions hold:
\begin{enumerate}
\item
if $\cK$ is acute, $\ccK$ is contained in the triangle formed by the segments
joining the midpoints of the sides of $\cK;$ also, $\ccK$ equals the
quadrangle $\cQ$ formed by the bisectors of the smallest and largest angles and 
the axes of the shortest and longest sides of $\cK;$
\item
if $\cK$ is obtuse, $\ccK$ is contained in the parallelogram whose vertices are
the midpoints mentioned in (i) and the vertex of the smallest angle in $\cK;$
also, $\ccK$ equals the polygon $\cP$ formed by the largest side of $\cK$ and
the bisectors and axes mentioned in (i); $\cP$ may be either a pentagon or a quadrangle.
\end{enumerate}
\end{lemma}
\begin{proof}
Observe that bisectors of angles and axes of sides are admissible axes of reflection.
If $\cK$ is acute, $C_\cK$ and $I_\cK$ fall in its interior and are the intersection of the axes
and bisectors, respectively. If $\cK$ is obtuse, $I_\cK$ still falls in the interior of $\cK,$
while $C_\cK$ is the midpoint of the largest side of $\cK$ and is no longer the
intersection of the axes. These remarks imply that $\ccK\subseteq\cQ$ in case (i) and $\ccK\subseteq\cP$
in case (ii); also $C_\cK, I_\cK\in\cQ\cap\ccK$ and $C_\cK, I_\cK\in\cP\cap\ccK.$ 
\par 
The segments 
specified in (i) are also admissible axes of reflection if $\cK$ is acute;
thus, the inclusion in the triangle mentioned in (i) easily follows.
If $\cK$ is obtuse, only the segment joining the midpoints of the smallest and intermediate side
is an axis of reflection. However, we can still claim that $\ccK$ is contained in the
parallelogram mentioned in (ii), since $C_\cK$ is now the midpoint of the largest side
from which one of the axes is issued: thus, $\ccK$ must stay below that axis.
\par
Now, if $\ccK$ were smaller than $\cQ$ (or $\cP$), then there would be an
axis of reflection that cuts off one of the vertices of $\cQ$ (or $\cP$) 
different from $C_\cK$ and $I_\cK$ (that always belong to $\ccK$). 
In any case, such an axis would violate the maximality that
axes of sides and bisectors of angles enjoy with respect of reflections.
\end{proof}

\par
We are now ready to prove Theorem \ref{th:triangle}.

\begin{proof} 
First of all, thanks to the inclusion mentioned in Lemma \ref{th:hearttriangle}, we get $|\ccK|<1/4\, |\cK|$ when $\cK$ is acute. 
Thus, we can restrict ourselves to the case of $\mathcal{K}$ obtuse.
\par 
Here, we refer to Figure \ref{fig3}. We observe that, by what we proved in Lemma \ref{th:hearttriangle}, 
$\ccK$ is always contained in the quadrangle $DEFG$ (when the angle in $B$ is much larger than $\pi/2,$ $DEFG$ and $\ccK$
coincide),
which is contained in the trapezoid $DELH.$ Thus, $|\ccK|\le |DELH|;$
hence it is enough to prove that, if the angle in $B$ increases, $|DELH|$ increases  and
both ratios $|\ccK|/|\cK|$ and $|DELH|/|\cK|$ tend to $3/8.$
\par
We proceed to compute $|\ccK|,$ when the angle in $B$ is large.  
We fix a base and a height of $\cK:$ as a base we choose the smallest side
and we suppose it has length $b;$ $h$ will denote the length of its corresponding height.
In this way, $|\cK|=b\,h/2.$

\begin{figure}
\label{fig3}
\includegraphics[scale=.4]{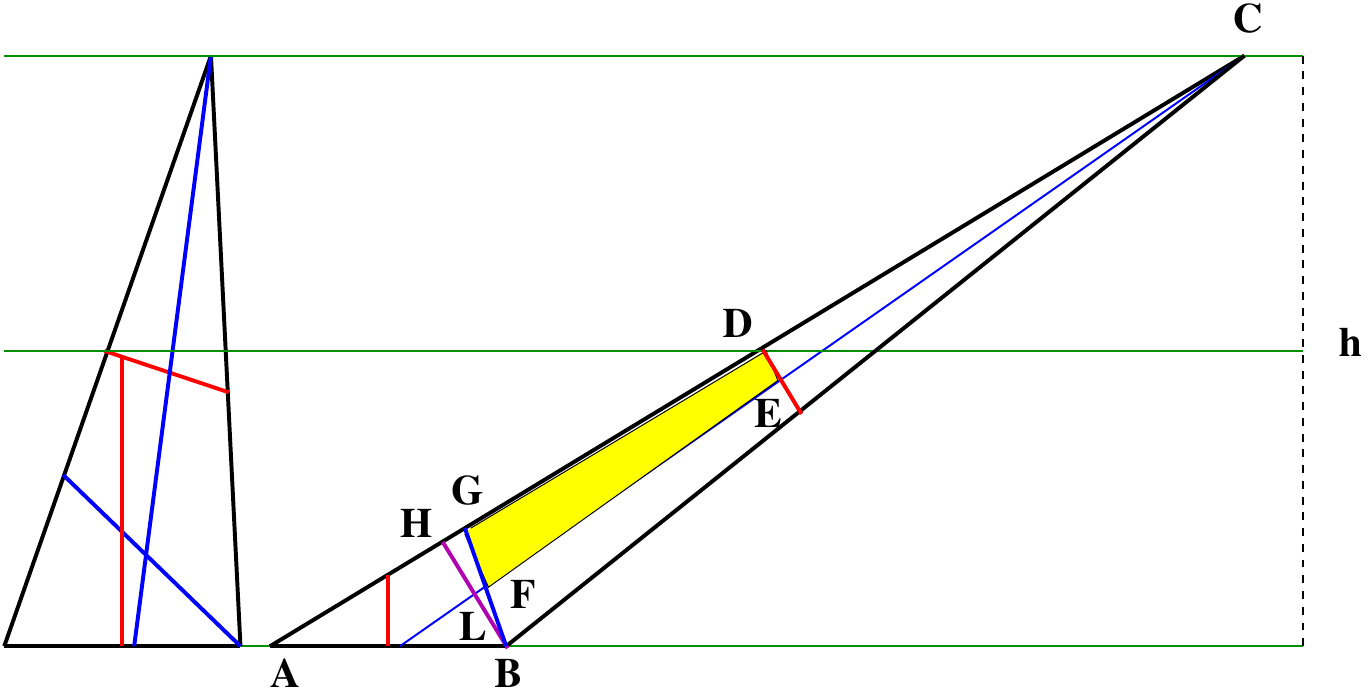}
\caption{The heart of a triangle.}
\end{figure}

\par
In Figure \ref{fig3}, the lines through the points $B$ and $G,$ and $C$ and $L$ 
bisect the angles in $B$ and $C,$ respectively. The line through $D$ and $E$
is the only axis that contributes to form $\ccK,$ that equals the quadrangle $DEFG;$ 
thus, $\ccK$ is obtained as
$$
\ccK=T_1\setminus (T_2\cup T_3),
$$
where the $T_i$'s are triangles: 
$$
T_1=CBG, \ T_2=CBF, \ T_3=CED.
$$
\par
We place the origin of cartesian axes in $A$ and set $B=(b,0);$
we also set $C=(t,h).$ Finally, we denote by $\al,$ $\be$ and $\ga$ the 
respective measures of the angles in $A,$ $B$ and $C.$ 
\par
Trigonometric formulas
imply that
\begin{eqnarray}
\label{triangles}
&&|T_1|=\frac12 \, [h^2+(t-b)^2]\,\frac{\sin(\be/2) \sin(\ga)}{\sin(\be/2+\ga)},\nonumber\\
&&|T_2|=\frac12\, [h^2+(t-b)^2]\, \frac{\sin(\be/2) \sin(\ga/2)}{\sin(\be/2+\ga/2)},\\
&&|T_3|=\frac18 \, [h^2+t^2]\,\tan(\ga/2)\nonumber,
\end{eqnarray}
where the angles $\be$ and $\ga$ are related by the {\it theorem of sines:}
$$
\frac{h}{\sqrt{h^2+t^2}\,\sqrt{h^2+(t-b)^2}}=\frac{\sin(\be)}{\sqrt{h^2+t^2}}=\frac{\sin(\ga)}{b}.
$$
The area of $DELH$ is readily computed as
\begin{equation}
\label{trapezio}
|DELH|=\frac12\,\frac{b^2\,h^2}{(h^2+t^2)\,\tan(\ga/2)}-\frac18\,(h^2+t^2)\,\tan(\ga/2).
\end{equation}
\par
Now, observe that this quantity increases with $t,$ since it is the composition of
two decreasing functions: $s\mapsto b^2\,h^2/(2s)-s/8$ and $t\mapsto (h^2+t^2)\,\tan(\ga/2).$
\par
As $t\to\infty,$ $|\cK|$ does not change, the angle $\ga$ vanishes
and the angle $\be$ tends to $\pi;$ moreover, we have that
$$
t\,\sin(\be)\to h, \ \ t^2\,\sin(\ga)\to bh \ \mbox{ as } \ t\to\infty.  
$$
Formulas \eqref{triangles} then yield:
$$
|T_1|\to \frac12\, bh=|\cK|, \ \ |T_2|\to\frac14\, bh=\frac12\,|\cK|, \ \ |T_3|\to\frac1{16}\, bh=\frac18\,|\cK|.
$$
Thus, since $|\ccK|=|T_1|-|T_2|-|T_3|,$ we have that $|\ccK|\to \frac38\,|\cK|;$ 
by \eqref{trapezio}, $|DELH|\to \frac38\,|\cK|$ as well. 
\par
The proof is complete.
\end{proof}

\begin{remark}
Thus, Theorem \ref{th:triangle} sheds some light on problem \eqref{shape}.  
In fact, observe that the maximizing sequence, once properly re-scaled, gives a maximizing sequence for the equivalent problem \eqref{diameter}, that precisely collapses to a one-dimensional object. 
\par
Numerical evidence based on the algorithm developed in
\cite{BMS} suggests that, the more $\cK$ is {\it round,} the more $\ccK$ is small
compared to $\cK.$ We conjecture that  
$$
\sup\left\{ \frac{|\ccK|}{|\cK|}: \cK\subset\RE^2 \mbox{ is a convex body} \right\}=\frac{3}{8},
$$
and the supremum is realized by a sequence of obtuse triangles.
\end{remark}

\begin{acknowledgement}
The authors wish to warmly thank Prof. J. O'Hara, who pointed out a gap in a preliminary version of the proof of Theorem \ref{teo:pmoment} and provided a copy of \cite{OH}.
The first author has been partially supported by the ERC Advanced Grant n. 226234.
The second author has been supported by GNAMPA-INdAM and the PRIN-MIUR grant
``Propriet\` a e metodi geometrici nelle equazioni alle derivate parziali, disuguaglianze di Sobolev e convessit\` a''. Both authors gratefully acknowledge the Banff International Research Station and its facilities, where this work started.
\end{acknowledgement}

\end{document}